\documentclass{article}
\usepackage{amsmath, amssymb}
\usepackage{graphicx,latexsym}

\renewcommand{\abstractname}{Abstract.}
\renewcommand\abstract{\hfil\break\topsep=0pt\partopsep=0pt\parsep=0pt\itemsep=0pt\relax
\trivlist\item[\hskip\labelsep
{\bfseries\abstractname}]\if!\abstractname!\hskip-\labelsep\fi}

\newcommand{\email}[1]{{(e-mail: #1)}}

\def\keywordname{{\bfseries Key words:}}
\def\keywords#1{\par\addvspace\baselineskip\noindent\keywordname\enspace
\ignorespaces#1}

\def\jelclassname{{{\bfseries JEL Classification:}} }
\def\jelclass#1{\par\addvspace\medskipamount\noindent\jelclassname\
\ignorespaces#1}

\def\subclassname{{\bfseries Mathematics Subject Classification (1991):} }
\def\subclass#1{\par\addvspace\medskipamount\noindent\subclassname\
\ignorespaces#1}

\def\title#1{\hfil\break\hfil\break
\hfil\break\par\addvspace\baselineskip\noindent
\ignorespaces{\LARGE\bf#1}\hfil\break}

\def\author#1{\par\addvspace\baselineskip\noindent
\ignorespaces{\large\bf#1}}

\def\institute#1{\par\addvspace\baselineskip\noindent
\ignorespaces{\small#1}\hfil\break}

\newtheorem{lemma}{Lemma}[section]

\newtheorem{remark}{Remark}

\newtheorem{theorem}{Theorem}[section]

\newtheorem{thm}{Theorem}[section]
\newtheorem{assumption}[thm]{Assumption}
\DeclareMathOperator*{\trace}{trace}
\DeclareMathOperator*{\argmin}{argmin}
\DeclareMathOperator*{\diag}{diag}
\DeclareMathOperator*{\sign}{sign}

\begin{document}

\title{Tracking of historical volatility}

\author{Goldentayer, L.}
\institute{Electrical Engineering Systems, Tel Aviv University,
69978 - Ramat Aviv, Tel Aviv, Israel}
\email{goldlev@hotmail.com}

\author{Klebaner F.}
\institute{School of Mathematical Sciences, Building 28M, Monash
University, Clayton Campus, Victoria 3800, Australia.}
\email{fima.klebaner@sci.monash.edu.au}

\author{Liptser, R.}
\institute{Electrical Engineering Systems, Tel Aviv University,
69978 - Ramat Aviv, Tel Aviv, Israel}
\email{liptser@eng.tau.ac.il}

\begin{abstract}
We propose an adaptive algorithm for tracking of historical
volatility. The algorithm is built under the assumption that the
historical volatility function belongs to the
Stone-Ibragimov-Khasminskii class of $k$ times differentiable
functions with bounded highest derivative and its subclass of
functions satisfying a differential inequalities. We construct an
estimator of the Kalman filter type and show optimality of the
estimator's convergence rate  to zero as sample size $n\to\infty$.
This estimator is in the framework of GARCH design, but
 a tuning procedure of its parameters
is faster than with traditional GARCH techniques.
\end{abstract}

\keywords{GARCH, Historical volatility, Volatility estimation,
On-line tracking estimator,  Global adaptation}

\jelclass{C13}

\subclass{60G35,60G51,62G05,62M20,91B70}

\newpage
\section{Introduction}
\label{sec-1}

\subsection{}
In the classical Black-Scholes model for financial markets, the
stock price $S(t)$ is modelled as a Geometric Brownian motion with
the diffusion coefficient ``$\sqrt{v }S(t)$'', where {\em
volatility} $v$ is assumed to be constant. This assumption is
convenient for the option price ``prediction''. Contrary to this
assumption, the traders treat the volatility as a parameter that
changes with time and whose future values have to be evaluated
(predicted) for a given period of interest. In this connection,
many researchers would rather interpret the volatility as a random
process, $v(t)$, and study, so-called, {\em stochastic} volatility
models. It is natural to verify how the volatility $v(t)$
 changes in time for
real stock prices and attempt to select a suitable stochastic
volatility model. Traditionally it is proposed to apply {\it
Generalized Autoregressive Conditional Heteroscedasticity}
tracking algorithms (shortly GARCH, see, e.g.
\cite{Torben}-\cite{Bol86}, \cite{Day}, \cite{Duan}, \cite{Ham94})
for tracking $v(t)$ from the stock prices. It is known from
\cite{Ham94} (p.109) that GARCH algorithm operates satisfactory
under relatively stable market conditions but fails when highly
unanticipated events that lead to a significant structural change
occur. Nevertheless, in many realistic settings, the simplest
GARCH(p,q), p,q=1,2, algorithms are adequate for tracking
volatilities even over long periods (see Bollerslev, Chou, and
Kroner \cite{BCK}, pages 10 and 22). The main difficulty in
implementation
 of GARCH comes from the multivariate minimization
procedure of  its parameters even for small values of p,q=1,2.

In this paper, we propose a new approach for tuning the GARCH
parameters. Our approach uses ideas from Nonparametric Statistics
combined with the Kalman-Bucy filter representation of a GARCH
model. This representation enables us to select a GARCH model with
only one parameter, that practically achieves tracking accuracy of
GARCH(1,1). Moreover, the Kalman-Bucy version of GARCH(p,q) allows
for a considerable simplification of the minimization procedure in
the GARCH parameters.

\subsection{}
Let  $S({t_i})$, $i=0,1,\ldots,n$ be the sample, with
$t_i-t_{i-1}\equiv :\triangle$, of asset prices $S(t)$, $0\le t\le
T$, from the Black-Scholes model (see \cite{Black1}, \cite{Black})
\begin{equation}\label{00}
    dS(t)=\mu(t)S(t)dt+\sqrt{v(t)}S(t)dB_t, \ 0\le t\le T,
\end{equation}
where $B_t$ is a Brownian motion, $S(0)$ is the initial stock
price, $\mu(t)$ and $v(t)$ are strictly positive deterministic
functions, respectively. Denote by
$$
X_i=\frac{1}{\triangle}\ln^2\Big(\frac{S(t_i)}{S(t_{i-1})}\Big).
$$
the observed heteroscedasticity. Since
$$
    \ln \Big(\frac{S(t_i)}{S(t_{i-1})}\Big)=
    \int_{t_{i-1}}^{t_i}0.5\big(2\mu(s)-v(s)\big)ds+\int_{t_{i-1}}^{t_i}
    \sqrt{v(s)}dB_s,
$$
By  It\^o's formula we find that
\begin{equation}\label{XiDef}
\begin{aligned}
X_i&=\frac{1}{\triangle}\Big(\int_{t_{i-1}}^{t_i}0.5\big(2\mu(s)-v(s)\big)ds
+\int_{t_{i-1}}^{t_i}\sqrt{v(s)}dB_s\Big)^2
\\
&=\frac{1}{\triangle}\Big(\int_{t_{i-1}}^{t_i}\sqrt{v(s)}dB_s\Big)^2
+\frac{1}{\triangle}\Big(\int_{t_{i-1}}^{t_i}0.5\big(2\mu(s)-v(s)\big)ds\Big)^2
\\
&+\frac{2}{\triangle}\int_{t_{i-1}}^{t_i}0.5\big(2\mu(s)-v(s)\big)ds
\int_{t_{i-1}}^{t_i}\sqrt{v(s)}dB_s.
\end{aligned}
\end{equation}
The parameter $\triangle$ is usually small (for example, if the
stock prices are measured once a day for three consecutive years,
then $\triangle \simeq 0.001$). For sufficiently small
$\triangle$, the dominating term in $X_i$ is
$\frac{1}{\triangle}(\int_{t_{i-1}}^{t_i}\sqrt{v(s)}dB_s)^2$, with
the mean
$$
v_{i-1}=\frac{1}{\triangle}\int_{t_{i-1}}^{t_i}v(s)ds
$$
 and, under some smoothness assumptions,
 the error $
\frac{1}{\triangle}\int_{t_{i-1}}^{t_i}(v(s)-v_{i-1})ds $ is
sufficiently small and can be ignored.

\medskip
Following Bollerslev \cite{Bol86} and Engle \cite{Eng82}, the
GARCH(p,q) provides estimates $\widehat{v}_i$ of $v_i$ by the
recursion
\begin{equation}\label{garch}
\widehat{v}_i = K+\sum_{j=1}^p g_j \widehat{v}_{i-j}+\sum_{m=1}^q
a_m X_{(i+1)-m},
\end{equation}
subject to some initial conditions, where parameters
$K,g_1,\ldots,g_p, a_1,\ldots, a_q$, as well as $p,q$, have to be
chosen with the help of minimizing the observed sum of squares
(here $n=\frac{T}{\triangle}$)
$$
S_n(K,g_1,\ldots,g_p, a_1,\ldots, a_q)
=\frac{1}{n}\sum_{i=1}^n(X_i-\widehat{v}_{i-1})^2.
$$

In contrast to \eqref{garch}, we propose an alternative tracking
algorithm borrowed, from Khasminskii - Liptser, \cite{KhL}, and
Goldentayer - Liptser \cite{GL}, with the univariate minimizing
parameter $\gamma$:
\begin{equation}\label{farsh}
\begin{aligned}
\widehat{v}_i&=\widehat{v}_{i-1}+
\frac{1}{n}\widehat{v}^{(1)}_{i-1}+
\frac{q_0(\gamma)}{n^{2(k+1)/(2k+3)}}
\big(X_i-\widehat{v}_{i-1}\big)
\\
\widehat{v}^{(j)}_i&=\widehat{v}^{(j)}_{i-1}+
\frac{1}{n}\widehat{v}^{(j+1)}_{i-1}+ \frac{q_j(\gamma)}{
n^{(2(k+1)-j)/(2k+3)}} \big(X_i-\widehat{v}_{i-1}\big)
\\
j&=1,\ldots,k-1
\\
\widehat{v}^{(k)}_i&=\widehat{v}^{(k)}_{i-1}+
\frac{q_k(\gamma)}{n^{(k+2)/(2k+3)}}
\big(X_i-\widehat{v}_{i-1}\big).
\end{aligned}
\end{equation}
subject to some initial conditions, where $q_i(\gamma)$,
$i=0,1,\ldots,k$ are some prescribed functions and $\gamma$ is
chosen  to minimize
$$
S_n(\gamma)=\frac{1}{n}\sum_{i=1}^n(X_i-\widehat{v}_{i-1})^2.
$$
A choice of the parameter $k$ is imposed by the smoothness of
$v(t)$: $k=0$ when  $v(t)$ is  Lipschitz continuous (with a global
Lipschitz constant) while other positive values of $k$ are used
when $v(t)$ has a bounded $k$-th derivative.

It should be noted that the estimator, given in \eqref{farsh}, is
also of GARCH type. In Section \ref{sec-3a}, we give modification
of \eqref{farsh} which is compatible with GARCH(p,q).

The proposed estimator admits a fast optimization procedure and
enables to avoid local minima quite easily. Moreover, its
GARCH(p,q) modifications are always stable and possess faster
minimization than classical GARCH.

\section{Description of estimator. Quality of estimation}
\label{sec-2}
\subsection{Assumptions on preliminaries}

We assume that $v(t)$ is a smooth function. If $v(t)$ is a paths
of random process we assume that this random process and the
Brownian motion $B(t)$ are independent.

We use Nonparametric Statistics ideas for estimating a smooth
function observed in the presence of white noise. Although we use
some adaptive techniques, our method is different to that of
Spokoiny and Mercurio \cite{Spok}, where the volatility is
approximated by a piecewise constant function.

\begin{assumption}\label{ksmooth}
The volatility $v(t)$ is a strictly positive bounded function and
belongs to the
Ibragimov - Khasminskii - Stone subclass of functions {\rm (}see,
{\rm \cite{IK80}, \cite{IK81}} and {\rm \cite{St82})} namely,  $k$ times
differentiable {\rm (}$k=0$ included{\rm )} with Lipschitz continuous $k$-th
derivative.
\end{assumption}
In the accordance with this assumption, there exists a positive
number $L$ such that for any $i$
\begin{equation}\label{2.1a}
|v(t_i)-v(t_{i-1})|\le L\frac{\triangle^{1+k}}{(1+k)!}.
\end{equation}

\begin{assumption}\label{small}
$\mu(t)$ is a positive and bounded function.
\end{assumption}
Set
\begin{equation}\label{mu}
\mu_{i-1}=\frac{1}{\triangle}\int_{t_{i-1}}^{t_i}\mu(s)ds.
\end{equation}

Introduce
\begin{equation*}\label{xi}
\xi_i=\frac{1}{\sqrt{v_{i-1}\triangle}}\int_{t_{i-1}}^{t_i}\sqrt{v(s)}dB_s,
\ i\ge 1
\end{equation*}
and notice that then $(\xi_i)_{i\ge1}$ forms an i.i.d. sequence of
$(0,1)$-Gaussian random variables. From \eqref{XiDef}, it follows
that
$$
X_i=0.25\big(2\mu_{i-1}-v_{i-1}\big)^2\triangle +2\sqrt{\triangle
v_{i-1}}\big(\mu_{i-1}-\frac{1}{2}v_{i-1}\big)\xi_i
+v_{i-1}\xi^2_i.
$$
Denote
\begin{equation*}\label{willbe}
\begin{aligned}
&\eta_i=\sqrt{\triangle v_{i-1}}\big(2\mu_{i-1}+v_{i-1}\big)\xi_i
+v_{i-1}(\xi^2_i-1),
\\
&\theta_i(\triangle)=0.25\triangle\big(2\mu_{i-1}+v_{i-1}\big)^2.
\end{aligned}
\end{equation*}
By Assumptions \ref{ksmooth} and  \ref{small}, $\theta_i(\triangle)=O(\triangle)$ and
$(\eta_i)_{i\ge 1}$ forms a sequence of zero mean uncorrelated
random variables with
\begin{equation}\label{eta}
E\eta^2_i=\triangle v_{i-1}\big(2\mu_{i-1}+v_{i-1}\big)^2
+2v^2_{i-1}=:\sigma^2_i
\end{equation}
with $\sigma^2_i$'s the strictly
positive and bounded numbers. Thus $X_i$ possesses the following
structure:
\begin{equation}\label{basic}
X_i=v_{i-1}+\eta_i+\theta_i(\triangle).
\end{equation}

\subsection{A reductive model for $\mathbf{X_i}$}

We replace \eqref{basic} by a simpler model:
\begin{equation}\label{simpl}
X_i=v_{i-1}+\eta_i
\end{equation}
For this model, it follows from  Ibragimov, Khasminskii
\cite{IK80}, \cite{IK81} (see also Stone \cite{St82}), that there
exists a kernel type estimate $\widehat{v}_i$ of $v_i$, generated
by $(X_i)_{1\le i\le n}$, such that for any $i$
\begin{equation}\label{raten}
E(v_i-\widehat{v}_i)^2\le O\big(n^{-2(1+k)/(2k+3)}\big).
\end{equation}
It is also known from Khasminskii  and Liptser \cite{KhL} that the
rate in $n$, given in \eqref{raten}, remains valid for the on-line
estimate obtained with the help of recurrent algorithm given below
\begin{equation}\label{1.4tri}
\begin{aligned}
\widehat{v}_i&=\widehat{v}_{i-1}+
\frac{1}{n}\widehat{v}^{(1)}_{i-1}+ \frac{q_0}{ n^{2(k+1)/(2k+3)}}
\big(X_i-\widehat{v}_{i-1}\big)
\\
\widehat{v}^{(j)}_i&=\widehat{v}^{(j)}_{i-1}+
\frac{1}{n}\widehat{v}^{(j+1)}_{i-1}+ \frac{q_j}{
n^{(2(k+1)-j)/(2k+3)}} \big(X_i-\widehat{v}_{i-1}\big)
\\
j&=1,\ldots,k-1
\\
\widehat{v}^{(k)}_i&=\widehat{v}^{(k)}_{i-1}+
\frac{q_k}{n^{(k+2)/(2k+3)}} \big(X_i-\widehat{v}_{i-1}\big).
\end{aligned}
\end{equation}
More exactly, the above-mentioned rate in $n$ is preserved out of
the boundary layer $i\ge O\big(n^{-1/(2k+3)} \log n\big)$,
resulting from uncertainty in the initial conditions for
\eqref{1.4tri}, provided that (see \cite{KhL})
\begin{assumption}\label{roots}
All roots of the characteristic polynomial
\begin{eqnarray}
p^k(\lambda,\mathfrak{q})=\lambda^{k+1}+q_0\lambda^k+q_1\lambda^{k-1}+\ldots+q_{k-1}\lambda+q_k
\label{1.5}
\end{eqnarray}
are different and have negative real parts.
\end{assumption}

\subsection{Adaptive estimator design}
\label{sec-2.3} Out of the above-mentioned boundary layer $i\ge
O\big(n^{-1/(2k+3)} \log n\big)$
estimates $(\widehat{v}_i)_{i\ge 1}$ obey the following property
(see, \eqref{raten}):
\begin{equation}\label{main0}
\varlimsup_{n\to\infty} \sup_{v_i}E\big(v_i- \widehat{v}_i\big)^2
n^{2(k+1)/(2k+3)}\le C(\mathfrak{q}),
\end{equation}
where the supremum is taken over all $v_i$'s satisfying Assumption
\ref{ksmooth}. The parameter $C(\mathfrak{q})$ depends on a filter
gain $\mathfrak{q}$, the vector with entries $q_0,q_1,\ldots,q_k$.
So, preserving the rate in $n$, the asymptotic estimation accuracy
depends on $\mathfrak{q}$ chosen in the framework of Assumption
\ref{roots}. It is clear that a direct minimization of
$C\big(\mathfrak{q}\big)$ in $\mathfrak{q}$ may contradict
Assumption \ref{roots}. Goldentayer and Liptser, \cite{GL},
proposed an approach, based on the Kalman-Bucy filtering theory,
for minimization of $C(\mathfrak{q})$ while preserving Assumption
\ref{roots}. However, for this approach the assumption
$E\eta^2_i\equiv \sigma^2$ was used. Although in the case
considered here $E\eta^2_i\not\equiv \sigma^2$, we
shall still follow this methodology. For a known $\sigma^2$ and a
free parameter $\gamma$, set
\begin{equation}\label{gamma/sigma}
\vartheta=\frac{\gamma}{\sigma}.
\end{equation}
This parameter  was introduced in  \cite{GL} and plays a crucial
role in creating the filtering gain $\mathfrak{q}$ entries:
\begin{equation*}\label{ewe}
\begin{aligned}
& q_0(\vartheta)=U_{00}\vartheta^{1/k+1}
\\
& q_1(\vartheta)=U_{01}\vartheta^{2/k+1}
\\
&..............................
\\
& q_k(\vartheta)=U_{0k}\vartheta^{k/k+1}
\\
& q_k(\vartheta)=U_{0k}\vartheta,
\end{aligned}
\end{equation*}
where $U_{0j}$'s are entries of the first column of $U$ the
positive definite matrix, which is the unique solution of the
algebraic Riccati equation ($^*$ is the transposition symbol) $
aU+Ua^*+B-UA^*AU=0 $ with matrices
\begin{equation*}\label{Aa}
A=
  \begin{pmatrix}
    1 & 0 & 0 & \ldots & 0
  \end{pmatrix},
\ a=\begin{pmatrix}
  0 & 1 & 0 & \cdots & 0 \\
  0 & 0 & 1 & \cdots & \vdots \\
  \vdots & \vdots &  & \ddots & 0 \\
  0 & 0 & 0 & & 1 \\
  0 & 0 & 0 & \cdots & 0 \\
\end{pmatrix},
B=
  \begin{pmatrix}
    1\\
    0\\
    \vdots
    \\
    0
  \end{pmatrix}
\end{equation*}
of sizes $1\times(1+k)$, $(1+k)\times(1+k$, $(1+k)\times 1$
respectively. It is known from \cite{GL} that
$$
\mbox{
\begin{tabular}{|c|c|c|c|c|c|}
\hline
  k & $U_{00}$ & $U_{01}$ & $U_{02}$ & $U_{03}$ & $U_{04}$ \\
\hline
  0 & 1 & NA & NA & NA & NA \\
  1 & $\sqrt{2}$ & 1 & NA & NA & NA \\
  2 & 2 & 2 & 1 & NA & NA \\
  3 & $\sqrt{4+\sqrt{8}}$ & $2+\sqrt{2}$ & $\sqrt{4+\sqrt{8}}$ & 1 & NA \\
  4 & $1+\sqrt{5}$ & $3+\sqrt{5}$ & $3+\sqrt{5}$ & $1+\sqrt{5}$ & 1 \\ \hline
\end{tabular}.
}
$$
So in \cite{GL}, we deal with the estimator
\begin{equation}\label{1.4var}
\begin{aligned}
\widehat{v}_i(\vartheta)&=\widehat{v}_{i-1}(\vartheta)+
\frac{1}{n}\widehat{v}^{(1)}_{i-1}(\vartheta)+
\frac{U_{00}\vartheta^{1/k+1}} {
n^{2(k+1)/(2k+3)}}\big(X_i-\widehat{v}_{i-1}(\vartheta)\big)
\\
\widehat{v}^{(j)}_i(\vartheta)&=\widehat{v}^{(j)}_{i-1}(\vartheta)+
\frac{1}{n}\widehat{v}^{(j+1)}_{i-1}(\vartheta)+
\frac{U_{0j}\vartheta^{(j+1)/k+1}} {n^{(2(k+1)-j)/(2k+3)}}
\big(X_i-\widehat{v}_{i-1}(\vartheta)\big)
\\
j&=1,\ldots,k-1
\\
\widehat{v}^{(k)}_i(\vartheta)&=\widehat{v}^{(k)}_{i-1}(\vartheta)+
\frac{U_{0k}\vartheta}{n^{(k+2)/(2k+3)}}
\big(X_i-\widehat{v}_{i-1}(\vartheta)\big).
\end{aligned}
\end{equation}

\subsection{Global Adaptation}
We propose to use the estimator \eqref{1.4var} for tracking
$v_i$'s when the $X_i$'s are defined in the accordance with
\eqref{basic}, i.e.
$$
X_i=v_{i-1}+\eta_i+\theta_i(\triangle).
$$
The univariate minimization with the help of $\vartheta$
guarantees Assumption \ref{roots}. However, since the variance of
the noise is not  constant and, moreover, unknown, an evaluation
of $C(\mathfrak{q}(\vartheta))$, as in \cite{GL}, would be
difficult. Therefore, we follow GARCH-technique adaptive method
(see, e.g., \cite{BB}) and evaluate $
V_n(\vartheta)=\frac{1}{n}\sum_{i=1}^n(v_{i-1}-\widehat{v}_{i-1}(\vartheta))^2
$ via $
S_{n}(\vartheta)=\frac{1}{n}\sum_{i=1}^n\big(X_i-\widehat{v}_{i-1}(\vartheta)\big)^2.
$ We show that for $\triangle$ sufficiently small with probability
close to one
$$
S_n(\vartheta')>S_n(\vartheta'') \Rightarrow
V_n(\vartheta')>V_n(\vartheta'').
$$
A crucial role in proving this implication plays the
above-mentioned asymptotical estimate $ E(v_i-\widehat{v}_i)^2\le
O\big(n^{-2(1+k)/(2k+3)}\big), \ n\to\infty, $ which is valid not
only for $\theta_i(\triangle)\equiv 0$ but also when
$\theta_i(\triangle)=O(\triangle)$ (see Lemma \ref{theLemma} in
Appendix A).

\begin{theorem}\label{theo-2.1}
For sufficiently large $n$ and any $\vartheta'\ne\vartheta''$ and
any $\varepsilon>0$
\begin{multline}\label{ox}
P\Big(\Big|[S_n(\vartheta')-S_n(\vartheta'')]-
[V_n(\vartheta')-V_n(\vartheta'')]\Big|
>\varepsilon\Big)
\\
\le\varepsilon^{-2} O\big(n^{-(4k+5)/(2k+3)}\big).
\end{multline}
\end{theorem}
\begin{proof}
Taking into account $\theta_i(\triangle)=O(\triangle)=O(n^{-1})$,
we find that
\begin{multline*}
S_n(\vartheta)=V_n(\vartheta)+\frac{1}{n}\sum_{i=1}^n\eta^2_i+\frac{2O(n^{-1})}{n}
\sum_{i=1}^n\big(v_{i-1}-\widehat{v}_{i-1}(\vartheta)\big)\\
+\frac{2}{n}\sum_{i=1}^n\big[\big(v_{i-1}-\widehat{v}_{i-1}(\vartheta)\big)
+ O(n^{-1})\big]\eta_i  + O(n^{-2})
\end{multline*}
and, therefore,
\begin{multline*}
\big[S_n(\vartheta')-S_n(\vartheta'')\big]-
\big[V_n(\vartheta')-V_n(\vartheta'')\big]
\\
=\frac{2O(n^{-1})}{n}\sum_{i=1}^n\big[\big(v_{i-1}-\widehat{v}_{i-1}(\vartheta')
\big)-\big(v_{i-1}-\widehat{v}_{i-1}(\vartheta'')\big)\big]
\\
+\frac{2}{n}\sum_{i=1}^n
\big[\big(\widehat{v}_{i-1}(\vartheta')-\widehat{v}_{i-1}(\vartheta'')\big)
+O(n^{-1})\big]\eta_i+O(n^{-2}).
\end{multline*}
For notational convenience, set
\begin{multline*}
r^2:=E\Bigg(\frac{2}{n}\sum_{i=1}^n
\big[\big(\widehat{v}_{i-1}(\vartheta')-\widehat{v}_{i-1}(\vartheta'')\big)
+O(n^{-1})\big]\eta_i
\\
+\frac{2O(n^{-1})}{n}\sum_{i=1}^n\big(\widehat{v}_{i-1}(\vartheta')
-\widetilde{v}_{i-1}(\vartheta'')\big) +O(n^{-2})\Bigg)^2.
\end{multline*}
The use of $\big(\sum_{\ell=1}^3r_1\big)^2\le
3\sum_{\ell=1}^3r^2_\ell$ provides
$$
\begin{aligned}
r^2&\le 3\Bigg(\frac{4}{n^2}\sum_{i=1}^n
E\big[\big(\widehat{v}_{i-1}(\vartheta')-\widehat{v}_{i-1}(\vartheta'')\big)
+O(n^{-1})\big]^2E\eta^2_i
\\
&\quad\quad +\frac{4O(n^{-2})}{n^4}E\Big(
\sum_{i=1}^n\big(\widehat{v}_{i-1}(\vartheta')-\widehat{v}_{i-1}(\vartheta'')\big)
\Big)^2+O(n^{-4})\Bigg) \equiv 3\sum_{\ell=1}^3r^2_\ell.
\end{aligned}
$$
Applying obvious estimates
\begin{equation*}
\begin{aligned}
&
E\big(\widehat{v}_{i-1}(\vartheta'')-\widehat{v}_{i-1}(\vartheta')\big)^2
\le 2E\big(v_{i-1}-\widehat{v}_{i-1}(\vartheta')\big)^2
+2E\big(v_{i-1}-\widehat{v}_{i-1}(\vartheta'')\big)^2
\\
&\quad\quad\quad\quad=
  \begin{cases}
    O(1), & i\le O(n^{-1/n^{2k+3}}\log n)
    \\
O\big(n^{-(2(k+1)/(2k+3)}\big), & i> O(n^{-1/n^{2k+3}}\log n),
  \end{cases}
\end{aligned}
\end{equation*}
we  get the following upper bounds for $ r^2_\ell, \ \ell=1,2,3$:

\medskip

\indent\indent\indent\indent
    $
r^2_1\le O(n^{-1})\Big(O(n^{-2})+\frac{1}{n}\sum_{i=1}^n
E\big(\widehat{v}_{i-1}(\vartheta')
-\widehat{v}_{i-1}(\vartheta'')\big)^2\Big), $

\indent\indent\indent\indent
 $ r^2_2=\frac{12O(n^{-2})}{n^4}E\Big(
\sum_{i=1}^n\big(\widehat{v}_{i-1}(\vartheta')-\widehat{v}_{i-1}(\vartheta'')\big)
\Big)^2 , $

\indent\indent\indent\indent  $r^2_3=O(n^{-2})\le
O\big(n^{-(4k+5)/(2k+3)}\big)$.

\medskip
\noindent Hence, with the help of  Chebyshev's inequality we find
that for sufficiently large $n$, any $\vartheta'\ne\vartheta''$
and any $\varepsilon>0$ the desired statement holds true.
\end{proof}

\begin{remark}
{\rm Theorem \ref{theo-2.1} enables a meaningful comparison
between estimators corresponding to various values of
 $\vartheta$. For notational convenience, the filter for $k=0$ with the best
parameter $\vartheta$ is called Filter 0.}
\end{remark}

\section{Filters controlled by multiple parameters}
\label{sec-3a}

We restrict ourselves by consideration of GARCH(1,1) and
GARCH(2,2) in the form of \eqref{1.4var}. To distinguish
these filters from classical GARCH's we denote them, by an analogy
with Filter 0, by Filter 1 and Filter 2, respectively. The
structure of these filters and  the motivation for their
applicability is given in Appendix \ref{App-B}. So, due to
\eqref{1.4var} and \eqref{1.5var}, we have

 {\bf Filter 1}
\begin{equation}\label{3.1c}
\widehat{v}_i=\widehat{v}_{i-1}
\Big(1-\frac{a_1}{n}\Big)+\frac{a_1 K}{n} +
\frac{\vartheta}{n^{2/3}} \big(X_i-\widehat{v}_{i-1}\big);
\end{equation}

{\bf Filter 2}
\begin{equation}\label{3.2c}
\begin{aligned}
\widehat{v}_i&=\widehat{v}_{i-1}+
\frac{1}{n}\widehat{v}^{(1)}_{i-1}+ \frac{\sqrt{2\vartheta}} {
n^{4/5}}\big(X_i-\widehat{v}_{i-1}\big)
\\
\widehat{v}^{(1)}_i&=\widehat{v}^{(1)}_{i-1}
\Big(1-\frac{a_1}{n}\Big)-\frac{a_2}{n}\widehat{v}_{i-1}+
\frac{a_2K}{n}
\\
&\quad + \frac{\vartheta}{n^{3/5}}
\big(X_i-\widehat{v}_{i-1}\big).
\end{aligned}
\end{equation}
It is assumed that $0<a_1, a_2 \ll n$ and $|K|\ll n$. The
estimates generated by Filters 1 and 2 possess the optimal rate in
$n\to\infty$, while for fixed $n$ the presence of additional
parameters $a_1, K$ and $a_1,a_2K$, respectively, enables slightly
to improve (about 10\%) the best value of
$$
S_n(\vartheta,K,
a_1,a_2)=\frac{1}{n}\sum_{i=1}\big(X_i-\widehat{v}_{i-1}\big)^2.
$$

The main adaptive parameter remains $\vartheta$. A contribution of
$a_1,K$ or $a_1,a_2,K$ is not  essential.
This fact enables to simplify
the tuning parameters procedure, particularly to avoid local
minima,
 in comparison with  the standard tuning procedure for classical
GARCH(1,1), GARCH(2,2) (see MATLAB GARCH Toolbox:
http://www.mathworks.com/access/helpdesk/help/toolbox/garch/garch.shtml

\section{Computer implementation and simulations}
\label{NumRes}

The volatility dynamics may differ widely between various types of
assets.  For example, the volatility changes for stocks and risky
assets are too fast and  the volatility values are relatively
high. The composite indexes and exchange rates characterized by
slow changes and smaller volatility values.
This remark points out the difficultly of finding the best filter
simultaneously for all assets.

In simulations, we compare the results of Filter 0, Filter 1 and
Filter 2, as well as the GARCH(1,1) and GARCH(2,2), provided by
MATLAB. Though Filter 1 and Filter 2 are equivalent to GARCH(1,1)
and GARCH(2,2) respectively,
the comparison of the numerical results  show some
advantage of Filters 1 and 2 due to different tuning procedures (see comment at the end of Section
\ref{sec-3a}).

\subsection{Tuning procedure for Filters 1 and 2}

The univariate minimization process required for Filter 0 is
straightforward. For Filter 1 and Filter 2 we used unconstrained
minimization, as given below.

\medskip

{\bf Filter 1} \mbox{}

1. Set $a_1=0$, $K=0$ and find $\vartheta^*=\argmin_\vartheta S_n(\vartheta,0,0)$.

2. Find $K^*=\frac{1}{n}\sum_{i=1}^nX_i$.

3. Find $a^*_1=\argmin_{a_1}S_n(\vartheta^*,K^*,a_1)$.

4. Local minimization of $S_n(\vartheta,K,a_1)$ in  vicinity of
$(\vartheta^*,K^*,a^*_1)$.

\medskip

{\bf Filter 2} \mbox{}

1. Set $a_1=0$, $a_2=0$, $K=0$ and find
$\vartheta^*=\argmin_\vartheta S_n(\vartheta,0,0,0)$.

2. Find $K^*=\frac{1}{n}\sum_{i=1}^nX_i$.

3. Find $(a^*_1,a^2_2)=\argmin_{a_1,a_2}S_n(\vartheta^*,K^*,a_1,a_2)$.

4. Local minimization of $S_n(\vartheta,K,a_1,a_2)$ in  vicinity
of $(\vartheta^*,K^*,a^*_1,a^*_2)$.

\medskip

The tuning procedures above consist in the univariate minimization
over $\vartheta$, the computation of $K$ and the
minimization over $a_1,a_2$. These steps are supposed to provide
some $(\vartheta^*,K^*,a^*_1,a^*_2)$ in the vicinity on the
minimum point, where the multidimensional minimization procedure is applied.

\begin{figure}[hbt]
\begin{center}
\includegraphics[angle=0,width=4.0in,height=3in]
{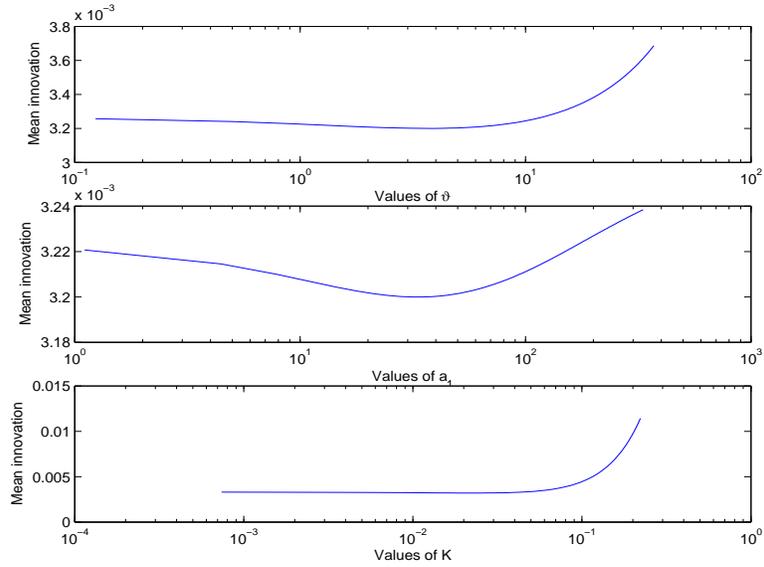}
\end{center}
\caption{Filter 1 coefficients behavior in the minima for IBM
stock.} \label{fig_c1}
\end{figure}

\begin{figure}[hbt]
\begin{center}
\includegraphics[angle=0,width=4.0in,height=3.5in]
{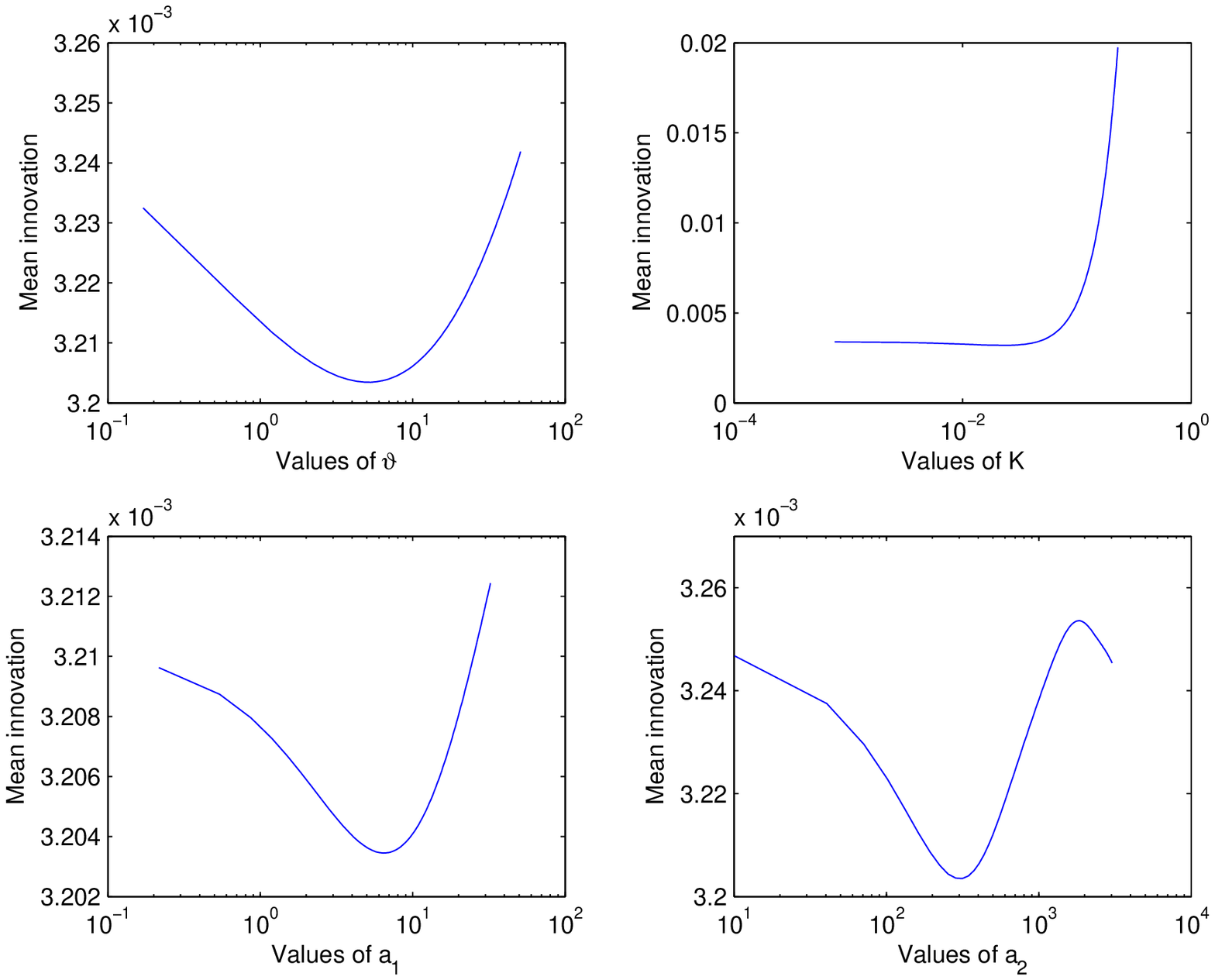}
\end{center}
\caption{Filter 2 coefficients behavior in the minima for IBM
stock.} \label{fig_c2}
\end{figure}

The simulation results demonstrate that in the vicinity of
$(\vartheta^*,K^*,a^*_1,a^*_2)$ the function
$S_n(\vartheta,K,a_1,a_2)$ behaves as a concave function and this
property is preserved in a wide range around of
$(\vartheta^*,K^*,a^*_1,a^*_2)$. Moreover, the minimum is not
sharp, so that the minimization procedure does not require high
'resolution'. The corresponding marginal projections of
$S_n(\vartheta,K,a_1,a_2)$ are given on Figures \ref{fig_c1} and
\ref{fig_c2}.


\subsection{Exchange rates}
The USD exchange rates we used for historical volatility
estimation were taken for the period between 01-Dec-01 and
18-Jan-04, e.g. $n=1466$.
Filter 1 and Filter 2 provide
estimation error similar to GARCH(1,1) and GARCH(2,2).
The tracking accuracy in
terms of $S_n(\vartheta,K,a_1,a_2)$ is given the Table \ref{Tab0}.

\begin{small}
\begin{table}[ht]  \centering    \caption{Average one-step prediction error for
 the exchange rates volatility}\label{Tab0}
 \begin{tabular}{|c|c|c|c|c|c|c|}
   \hline   \multicolumn{2}{|c|}{Currency} & \multicolumn{5}{|c|}{Filter type} \\
       \cline{1-7} From  & To & Garch(1,1) & Garch(2,2) & Filter 0 & Filter 1 & Filter 2\\  \hline
AUD & \$  & 9.096e-006 & 9.090e-006 & 9.117e-006 & 9.092e-006 & 9.092e-006\\
EUR & \$  & 5.857e-006 & 5.848e-006 & 5.869e-006 & 5.856e-006 & 5.856e-006\\
NIS & \$  & 1.830e-006 & 1.827e-006 & 1.839e-006 & 1.827e-006 & 1.826e-006\\
RUB & \$  & 5.026e-007 & 4.915e-007 & 5.026e-007 & 4.935e-007 & 4.880e-007\\
YEN & \$  & 5.388e-006 & 5.372e-006 & 5.389e-006 & 5.376e-006 & 5.367e-006\\
 \hline    \end{tabular}     \end{table}
\end{small}

\begin{figure}[hbt]
\begin{center}
\includegraphics[angle=0,width=4.3in,height=4.3in]
{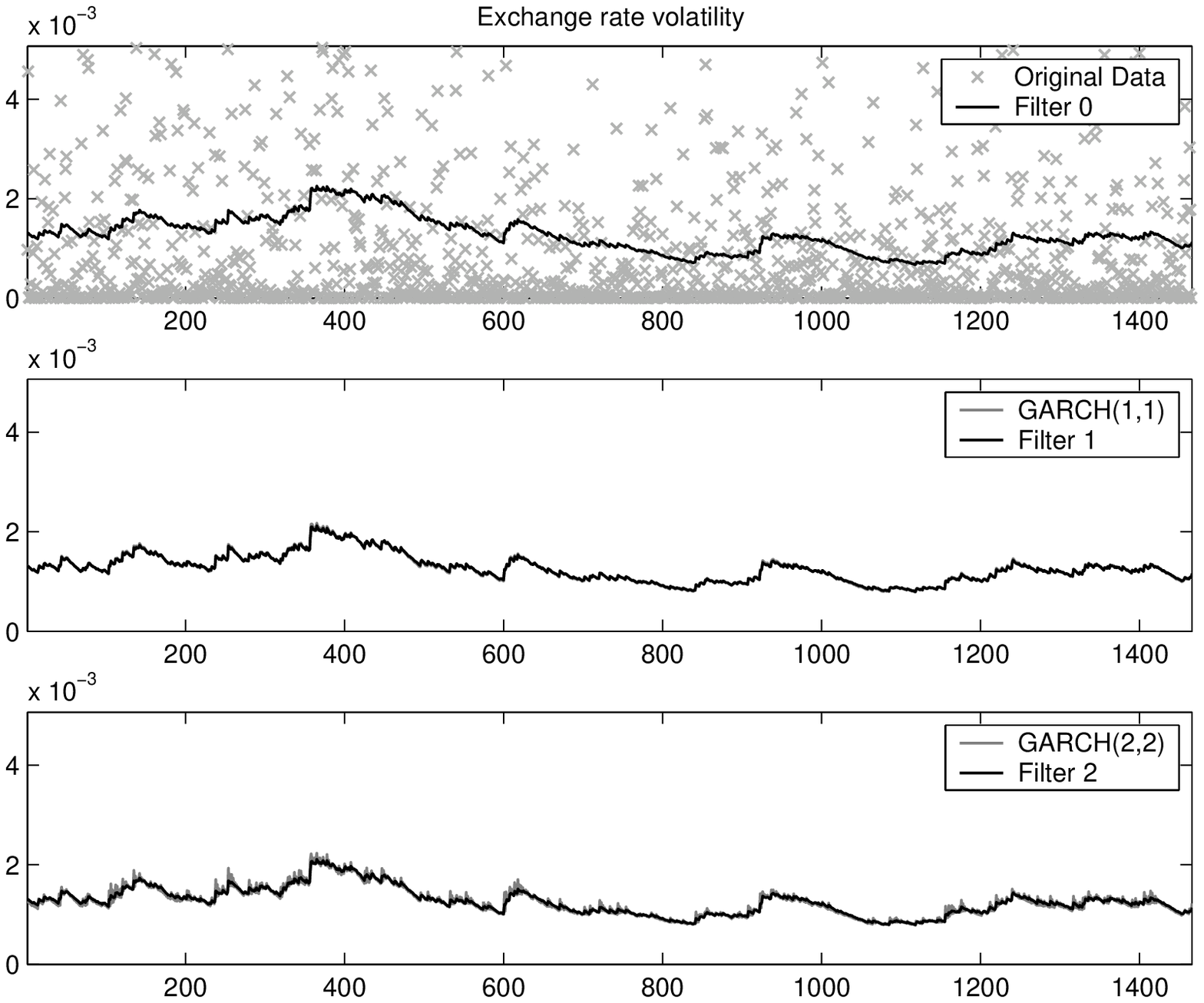}
\end{center}
\caption{Historical volatility estimation for exchange rates of
USD vs EUR.} \label{fig0}
\end{figure}

\subsection{Stocks}
We considered some stocks of large computer manufacturers and toy
and entertainment
  companies. The information (adjusted close
  prices) for the period: 24-Feb-99 to 28-Oct-03,
  $n=1176$, was collected using Yahoo.
The numerical results for $S_n(\vartheta,K,a_1,a_2)$,
corresponding the best tuning parameters, show serious differences
in filters behavior (see Table \ref{Tab1}).
A different quality of Filters 1, 2 and GARCH(1,1), (2,2) is
provided by different tuning procedures. Filter 1 provides the
best quality.
\begin{small}
 \begin{table}[h]  \centering    \caption{Average one-step prediction error for  the filter}\label{Tab1}
 \begin{tabular}{|l|c|c|c|c|c|}
   \hline   Asset  & \multicolumn{5}{|c|}{Filter type} \\     \cline{2-6} name  & Garch(1,1) & Garch(2,2) & Filter 0 & Filter 1 & Filter 2\\  \hline
DIS &  4.309e-003 & 4.314e-003 & 4.329e-003 & 4.284e-003 & 4.285e-003\\
HPQ &  3.741e-002 & 3.771e-002 & 3.615e-002 & 3.608e-002 & 3.608e-002\\
IBM &  3.229e-003 & 3.228e-003 & 3.217e-003 & 3.199e-003 & 3.200e-003\\
INTC & 1.232e-002 & 1.230e-002 & 1.235e-002 & 1.223e-002 & 1.232e-002\\
MAT &  1.899e-002 & 1.879e-002 & 1.850e-002 & 1.817e-002 & 1.847e-002\\
SUN & 5.144e-004 & 5.131e-004 & 5.143e-004 & 5.135e-004 & 5.138e-004\\
TOY & 6.134e-003 & 6.128e-003 & 6.158e-003 & 6.098e-003 & 6.082e-003\\
 \hline    \end{tabular}     \end{table}
\end{small}

\begin{figure}[hbt]
\begin{center}
\includegraphics[angle=0,width=4.3in,height=4.3in]
{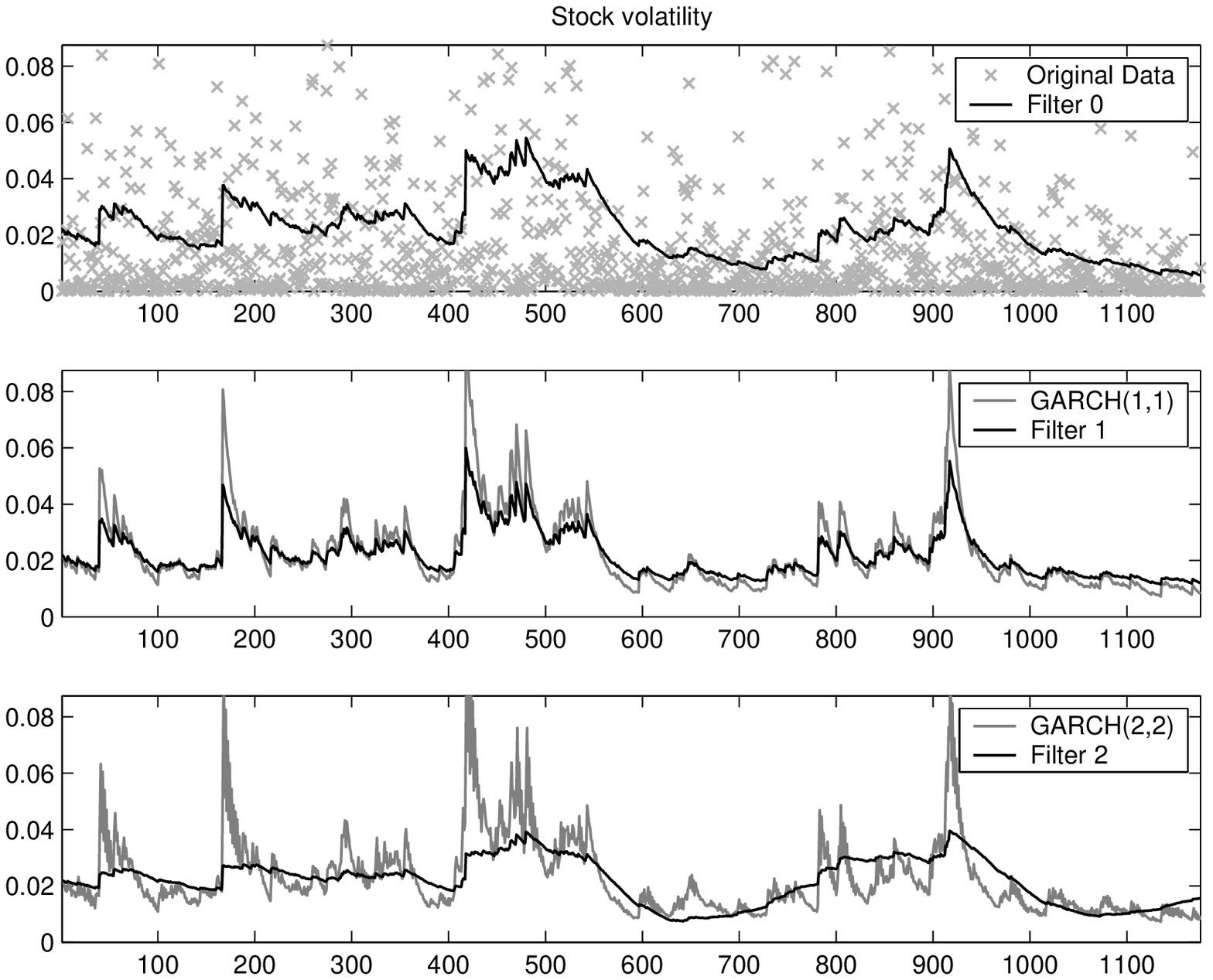}
\end{center}
\caption{Historical volatility estimation for IBM stock.}
\label{fig1}
\end{figure}

\subsection{Discussion of the numerical results}
The univariate minimization of
Filter 0 proves empirically to be very fast. The tuned parameter
$\vartheta^*$ for Filter 0 gives a hint for Filters 1 and  2
tuning procedures. The multivariate designs of Filters 1 and
  2 provide slightly better tracking accuracy than
GARCH(1,1) and (2,2), respectively, especially for stock options.
We attribute this effect to difficulties in tuning procedure of
the filter parameters, especially for GARCH(2,2), which leads to
local minima more often than tuning of for Filter 1 and 2.

\section{Conclusions}
The volatility estimation designs introduced in this paper are
similar to the widely used GARCH algorithms. The presented designs
allow simple adaptation and filter structure preserving
estimation accuracy similar to GARCH. The filter structure  of the presented design
 allows performance accuracy evaluation and enforces the stability of
the estimator.

\appendix
\section{Auxiliary Lemma}
\label{App-A}
\begin{lemma}\label{theLemma}
The rate in $n$ given in \eqref{raten} is preserved in the
presence of $\theta_i(\triangle)\not\equiv 0$ in \eqref{raten}.
\end{lemma}

\begin{proof}
Let $\widehat{v}_i$ and $\widetilde{u}_i$ be two estimates created
by \eqref{1.4var} with and without nuisance parameter
respectively. We prove below that
\begin{equation}\label{main**}
E\big(\widehat{v}_i-\widetilde{u}_i\big)^2\le O\big(n^{-2}\big).
\end{equation}
Since
$$ \frac{n^{2(k+1)/(2k+3)}}{n^2}\to 0, \ n\to\infty, $$
 the nuisance parameter $\theta_i(\triangle)$ does not change the
 rate \eqref{raten}
in $n\to\infty$.

\medskip
For the notational convenience write $\theta_i$, instead of
$\theta_i(\triangle)$, and set
$\delta^{(0)}_i=\widehat{v}_i-\widetilde{u}_i$,
$\delta^{(j)}_i=\widehat{v}^{(j)}_i-\widetilde{u}^{(j)}_i$,
$j=1,\ldots,k$. Then, we have
\begin{equation}\label{1.45}
\begin{aligned}
\delta_i&=\delta_{i-1}+ \frac{1}{n}\delta^{(1)}_{i-1}+
\frac{U_{00}\vartheta^{1/k+1}}{ n^{2(k+1)/(2k+3)}}
\big(\theta_i-\delta_{i-1}\big)
\\
\delta^{(j)}_i&=\delta^{(j)}_{i-1}+
\frac{1}{n}\delta^{(j+1)}_{i-1}+
\frac{U_{0j}\vartheta^{(j+1)/k+1}}{ n^{(2(k+1)-j)/(2k+3)}}
\big(\theta_i-\delta_{i-1}\big)
\\
j&=1,\ldots,k-1
\\
\delta^{(k)}_i&=\delta^{(k)}_{i-1}+
\frac{U_{0k}\vartheta}{n^{(k+2)/(2k+2)}}
\big(\theta_i-\delta_{i-1}\big)
\end{aligned}
\end{equation}
subject to the initial conditions $\delta(0)=0, \
\delta^{(j)}(0)=0, \  j=1,\ldots,k$.

Set
$$
\mathfrak{q_n}=\left(\begin{array}{c}
U_{00}\vartheta^{1/(1+k)}n^{-2(1+k)/(2k+3)}
\\
U_{01}\vartheta^{2/(1+k)}n^{-(2(1+k)-1)/((2k+3)}
\\
\vdots
\\
U_{0k}\vartheta n^{-(2(1+k)-k)/(2k+3)}
\end{array}
\right), \quad \mathfrak{q}=\left(\begin{array}{c}
U_{00}\vartheta^{1/(1+k)}
\\
U_{01}\vartheta^{2/(1+k)}
\\
\vdots
\\
U_{0k}\vartheta
\end{array}
\right)
$$
and $ F_i= \left(  \begin{array}{c}
    \delta_i
    \\
    \delta^{(1)}_i
    \\
    \vdots
    \\
    \delta^{(k)}_i
  \end{array}
  \right)
$ and recall that matrices $a$ and $A$ are defined in \eqref{Aa}.
We rewrite \eqref{1.45} to the vector-matrix form
$$
F_i=F_{i-1}+\frac{1}{n}aF_{i-1}+\mathfrak{q}_n\theta_i-
\mathfrak{q}_nAF_{i-1}
$$
where $F_0=0$. Set $G_i=C_nF_i$, where $C_n$ is the diagonal
$(1+k)\times(1+k)$-matrix:
\begin{small}
$$
C_n=\left(
\begin{array}{ccccccc}
    n^{(1+k)/(2k+3)} & 0 & \ldots & 0 & 0
    \\
    0 & n^{k/(2k+3)} & \ldots & 0 & 0
    \\
    \vdots & \vdots & \vdots & \vdots & \vdots
    \\
    0 & 0 & \ldots & n^{2/(2k+3)} & 0
    \\
    0 & 0 & \ldots & 0 & n^{1/(2k+3)}
    \end{array}
    \right).
$$
\end{small}
Then, $G_0=0$ and
\begin{equation}\label{3.5m}
G_i=G_{i-1}+\frac{1}{n}C_naF_{i-1}+C_n\mathfrak{q}_n\theta_i-
C_n\mathfrak{q}_nAF_{i-1}.
\end{equation}
By directly verifying identities $ C_na=n^{\frac{1}{2k+3}}aC_n$, \
$ C_n\mathfrak{q}_n=n^{-(1+k)/(2k+3)}\mathfrak{q}, $ we have
\begin{equation}\label{ustal}
\begin{aligned}
& \frac{1}{n}C_naF_{i-1}=n^{-2(1+k)/(2k+3)}aG_{i-1}
\\
& C_n\mathfrak{q}_n\theta_i=n^{-(1+k)/(2k+3)}\mathfrak{q}\theta_i.
\end{aligned}
\end{equation}
The structure of matrix $A$ provides $ n^{(1+k)/(2k+1)}A=AC_n. $
Hence and from $ C_n\mathfrak{q}_n=n^{-(1+k)/(2k+3)}\mathfrak{q},
$ it holds
\begin{equation}\label{3.7}
C_n\mathfrak{q}_nAF_{i-1}=n^{-2(1+k)/(2k+3)}\mathfrak{q}AG_{i-1}.
\end{equation}
Gathering now \eqref{3.5m}, \eqref{ustal}, \eqref{3.7}, we find
the recurrent equation for $G_i$'s:
$$
G_i=G_{i-1}+n^{-2(1+k)/(2k+3)}\big(a-\mathfrak{q}A\big)G_{i-1}
+n^{-(1+k)/(2k+3)}\mathfrak{q}\theta_i.
$$
With the matrix $ D_n=I+n^{-2(1+k)/(2k+3)}\big(a-\mathfrak{q}A), $
this recurrent equation is transformed into $ G_i=D_nG_{i-1}
+n^{-(1+k)/(2k+3)}\mathfrak{q}\theta_i. $ Hence, due to $G_0=0$,
we have $
G_i=\sum_{p=1}^iD^{i+1-p}_nn^{-(1+k)/(2k+3)}\mathfrak{q}\theta_p.
$ The latter and $|\theta_i|=O(n^{-1})$ provide
\[
\|G_i\|\le O(n^{-1})n^{-(1+k)/(2k+3)}\sum_{p=0}^\infty\|D^p_n\|
=O(n^{-(3k+4)/(2k+3)}\sum_{p=0}^\infty\|D^p_n\|.
\]
On the other hand, by Statement 2 in \cite{KhL}, for some positive
constants $c_\circ$ and $C$ and any $p$ the following estimate for
$\|D^p_n\|$ is valid:
\[
\|D^p_n\|\le C\exp\big(-c_\circ n^{-2p(1+k)/(2k+3)}\big).
\]
Consequently,
$$
\begin{aligned}
\|G_i\|&\le O\big(n^{-(3k+4)/(2k+3)}\big) \Big(1-e^{-c_\circ
n^{-2(1+k)/(2k+3)}}\Big)^{-1}
\\
&=O\big(n^{-(k+2)/(2k+3)}\big).
\end{aligned}
$$
Finally, by $F_i=C^{-1}_nG_i$, we find that
$$
|\delta^{(j)}_i|\le n^{-(1+k-j)/(2k+3)}\|G_i\|\le
O\big(n^{-1+j/(2k+3)}, \ j=0,1,\ldots,k.
$$

Hence, $ \delta^{(0)}_i\le O\big(n^{-2}\big) $ and \eqref{main**}
holds true.
\end{proof}

\section{GARCH in the form of (\ref{1.4var})}
\label{App-B}

The filter of \eqref{farsh} type was proposed in \cite{KhL} for
tracking of functions from the Ibragimov - Khasminskii - Stone
(IKS) class (see, \cite{IK80}, \cite{IK81} and \cite{St82}). A
further implementation of this filter compatible with
IKS(k)-subclasses, $k=0,1,\ldots$,  was developed in \cite{GL}
({\it function $f\in{\rm (IKS)(k)}$ if it is $k$-times
differentiable having Lipschitz continuous $k$-th derivative}).

An analysis of GARCH filter structure enables us to claim that
GARCH filter is compatible even  with further subclass of
(IKS)(k)'s. For fixed $k$, function $f$ from the subclass of
(IKS)(k) satisfies a differential inequality (see, Goldenshluger
and Nemirovski, \cite{GN}):

$$
\big|f^{(k)}(t)+a_1f^{(k-1)}(t)+a_2f^{(k-2)}(t)+\cdots+a_{k-1}f^{(1)}+a_kf(t)+a_k\big|
\le L,
$$
where $a_1,a_2,\ldots,a_k$ are such that the roots of polynomial
$$
P(x)=x^k+a_1x^{k-1}+a_2x^{k-2}+\cdots+a_k
$$
have negative real parts.

By an analogy to Chow, Khasminskii and Liptser, \cite{CKL}, we
propose the following estimator for $f(t)$ and its derivatives
$f^{(j)}(t)$, $j=1,\ldots,k$ via the observations of the process
$X_t$ with $X_0=0$ and $ dX_t=f(t)dt+\varepsilon dW_t: $
\begin{equation}\label{B.1b}
\begin{aligned}
d\widehat{f}(t)&=\widehat{f}^{(1)}(t)dt+
{\frac{q_0}{\varepsilon^{2/(2k+3)}}}
 \big(dX_t-\widehat{f}(t)dt\big),
\\
d\widehat{f}^{(j)}(t)&=\widehat{f}^{(j+1)}(t)dt+ {\frac{q_{j}}{
\varepsilon^{2j/(2k+3)}}} \big(dX_t-\widehat{f}(t)dt\big), \
j=1,...,k-1,
\\
d\widehat{f}^{(k)}(t)&=-\big(a_1\widehat{f}^{(k-1)}(t)+a_2\widehat{f}^{(k-2)}(t)
+\cdots+a_{k-1}\widehat{f}^{(1)}(t)+a_k\widehat{f}(t)
\\
&\quad +a_kK\big)dt +{\frac{q_{k}}{\varepsilon^{2k/(2k+3)}}}
\big(dX_t-\widehat{f}(t)dt\big).
\end{aligned}
\end{equation}
Now we show that
\begin{equation}\label{B.2b}
\begin{aligned}
E\big(f(t)-\widehat{f}(t)\big)^2&\asymp \varepsilon^{4k/(2k+3)}
\\
E\big(f^{(j)}(t)-\widehat{f}^{(j)}(t)\big)^2&\asymp
\varepsilon^{4(k-j)/(2k+3)}, \ j=1,...,k.
\end{aligned}
\end{equation}
Set $\triangle(t)=f(t)-\widehat{f}(t)$ and
$\triangle^{(j)}(t)=f^{(j)}(t)- \widehat{f}^{(j)}(t)$. From
\eqref{B.1b} and
$$
\begin{aligned}
\dot{f}(t)&=f^{(1)}(t)
\\
\dot{f}^{(j)}(t)&=f^{(j+1)}(t), \quad j=1,...,k-1,
\\
\dot{f}^{(k)}(t)&=-\big(a_1f^{(k-1)}(t)+a_2f^{(k-2)}(t)
+\cdots+a_{k-1}f^{(1)}(t)+a_kf(t)+a_kK\big),
\end{aligned}
$$
we derive
$$
\begin{aligned}
d\triangle(t)&=\triangle^{(1)}(t)dt-
{\frac{q_0}{\varepsilon^{2/(2k+3)}}}
 \big(\varepsilon dW_t+\triangle(t)dt\big),
\\
d\triangle^{(j)}(t)&=\triangle^{(j+1)}(t)dt- {\frac{q_{j}}{
\varepsilon^{2j/(2k+3)}}} \big(\varepsilon
dW_t+\triangle(t)dt\big), \ j=1,...,k-1,
\\
d\triangle^{(k)}(t)&=-\Big(\sum_{\ell=1}^{k-1}a_\ell\triangle^{(k-\ell)}(t)+
a_k\triangle(t)+u(t)\Big)dt
\\
&\quad -{\frac{q_{k}}{\varepsilon^{2k/(2k+3)}}} \big(\varepsilon
dW_t+\triangle(t)dt\big).
\end{aligned}
$$
Following \cite{CKL}, we introduce
\begin{equation}\label{B.3f}
\delta(t)={\frac{\triangle(t\varepsilon^{2/(2k+3)})}{\varepsilon^
{2(k+1)/(2k+3)}}}, \ \
\delta^{(j)}(t)={\frac{\triangle^{(j)}(t\varepsilon^{2/(2k+3)}) }
{\varepsilon^{2(k+1-j)/(2k+3)}}}, \ j=1,\ldots,k
\end{equation}
and notice that
\begin{equation}\label{B.4a}
\begin{aligned}
d\delta(t)&=[\delta^{(1)}(t)-q_0\delta(t)]dt-q_0dW^\varepsilon_t,
\\
d\delta^{(j)}(t)&=[\delta^{(j+1)}(t)-q_j\delta(t)]dt
-q_0dW^\varepsilon_t, \ j=1,...,k-1,
\\
d\delta^{(k)}(t)&=-
\Big(\sum_{\ell=1}^{k-1}\varepsilon^{2\ell/(2k+3)}a_\ell\delta^{(k-\ell)}(t)+
\varepsilon^{2k/(2k+3)}a_k\delta(t)\Big)dt
\\
&\quad
-u(t\varepsilon^{2/(2k+3)})dt-q_k\delta(t)dt-q_{k}dW^\varepsilon_t,
\end{aligned}
\end{equation}
where
$W^\varepsilon_t=\frac{1}{\varepsilon^{1/(2k+3)}}W_{t\varepsilon^{2/(2k+3)}}$
is the standard Wiener process.

Set
$$
D(t)=
\begin{pmatrix}
  \delta(t)
  \\
  \delta^{(1)}(t)
  \\
  \vdots
  \\
  \delta^{(k)}
\end{pmatrix},
\quad U(t)=\begin{pmatrix}
  0 \\
  \vdots \\
  0 \\
  u(t\varepsilon^{2/(2k+3)}
\end{pmatrix}
\quad Q=
\begin{pmatrix}
  q_0
  \\
  q_1
  \\
  \vdots
  \\
  q_k
\end{pmatrix}
$$
and introduce matrices $a_\varepsilon$ and $A$ of sizes
$(k+1)\times (k+1)$ and $1\times {k+1}$, respectively,
$$
a_\varepsilon =\begin{pmatrix}
0 & 1 & 0 & . & . & . & 0 & 0 &  \\
0 & 0 & 1 & . & . & . & 0 & 0 &  \\
. & . & . & . & . & . & . & . &  \\
. & . & . & . & . & . & . & . &  \\
. & . & . & . & . & . & . & . &  \\
0 & 0 & 0 & . & . & . & 0 & 1 &  \\
-\varepsilon^{2k/(2k+3)}a_k & . & . & . & . & . & . &
-\varepsilon^{2/(2k+3)}a_1 &
\end{pmatrix}
$$
and $ A=\begin{pmatrix}
  1 & 0 & \ldots & 0.
\end{pmatrix}.
$ With the help of the introduced matrices, we rewrite
\eqref{B.4a} into the vector-matrix form:
$$
dD(t)=(a_\varepsilon-QA)D(t)dt+U(t)dt-QdW^\varepsilon_t.
$$
It is known from \cite{CKL} that the vector $Q$ may be chosen such
that the eigenvalues of the matrix $a_0-QA$ have negative real
parts. This property is preserved for $a_\varepsilon-QA$, at least
for sufficiently small $\varepsilon$. Henceforth, we assume this
property for $a_\varepsilon-QA$ too. Consider now the Lyapunov
equation (here $\mathrm{I}$ is a unit matrix)
\begin{equation}\label{B.5a}
(a_\varepsilon-QA)P_\varepsilon+P_\varepsilon(a_\varepsilon-QA)^*+\mathrm{I}=0.
\end{equation}
It is clear that $P_\varepsilon$ is the unique positive definite
matrix and $\lim_{\varepsilon\to 0}P_\varepsilon=P_0$, where $P_0$
is the unique matrix solving the Lyapunov equation
$$
(a_0-QA)P_0+P_0(a_0-QA)^*+\mathrm{I}=0.
$$
Denote $\|D(t)\|^2_{P_\varepsilon}(:=\big\langle
D^*(t)P_\varepsilon D(t)\big\rangle)$. With the help of the It\^o
formula and \eqref{B.5a} we find that
$$
\begin{aligned}
d\|D(t)\|^2_{P_\varepsilon}&=\Big(-\|D(t)\|^2+2\big\langle
D(t),P_\varepsilon U(t)\big\rangle+ \langle Q,P_\varepsilon
Q\rangle\Big)dt
\\
&\quad +2\big\langle D(t),P\varepsilon QdW_t\big\rangle.
\end{aligned}
$$
Therefore the function $V(t)=E\|D(t)\|^2_{P_\varepsilon}$ is
differentiable and
$$
\dot{V}(t)=E\Big(-\|D(t)\|^2+2\big\langle D(t),P_\varepsilon
U(t)\big\rangle+ \big\langle Q,P_\varepsilon Q\big\rangle\Big).
$$
It is obvious that for sufficiently small $\varepsilon$ positive
constant $c_1,c_2,c_3$, can be found such that $E\|D(t)\|^2\ge
c_1V(t)$, $ 2E\big\langle D(t),P_\varepsilon U(t)\big\rangle\le
c_2\sqrt{V(t)}$, $ \big\langle Q,P_\varepsilon Q\big\rangle\le c_3. $
Hence, $ \dot{V}(t)\le -c_1V(t)+c_2\sqrt{V(t)}+c_3$. The use of the inequality
$\sqrt{x}\le\alpha^{-1}+\alpha x$, $\alpha\ge 0, x\ge 0$, with $\alpha=\frac{c_1}
{2c_2}$, provides
$$
\dot{V}(t)\le -0.5c_1V(t)+\frac{2c^2_2}{c_1}+c_3.
$$
Hence, for any $t\ge 0$ we have
$
V(t)\le V(0)+\frac{2c^2_2+c_1c_3}{0.5c^2_1}.
$

Obviously, this property and \eqref{B.3f} provide \eqref{B.2b}.

\medskip
For the discrete time setting with $\sqrt{t_i-t_{i-1}}\equiv\varepsilon$,
adapted to the framework of \cite{KhL}, we find that
$$
\begin{aligned}
\widehat{f}_i&=\widehat{f}_{i-1}+
\frac{1}{n}\widehat{v}^{(1)}_{i-1}+ \frac{q_0}
{ n^{2(k+1)/(2k+3)}}\big(X_i-\widehat{f}_{i-1}\big)
\\
\widehat{f}^{(j)}_i&=\widehat{f}^{(j)}_{i-1}+
\frac{1}{n}\widehat{f}^{(j+1)}_{i-1}+
\frac{q_1} {n^{(2(k+1)-j)/(2k+3)}}
\big(X_i-\widehat{f}_{i-1}\big)
\\
j&=1,\ldots,k-1
\\
\widehat{f}^{(k)}_i&=\widehat{f}^{(k)}_{i-1}
\Big(1-\frac{a_1}{n}\Big)-\frac{1}{n}
\Big(\sum_{\ell=2}^{k-1}a_\ell\widehat{f}^{(k-\ell)}_{i-1}+
a_k\widehat{f}_{i-1}+a_kK\Big)
\\
&\quad + \frac{q_k}{n^{(k+2)/(2k+3)}}
\big(X_i-\widehat{f}_{i-1}\big).
\end{aligned}
$$

\medskip
Finally, in the framework of this paper, we have
\begin{equation}\label{1.5var}
\begin{aligned}
\widehat{v}_i&=\widehat{v}_{i-1}+
\frac{1}{n}\widehat{v}^{(1)}_{i-1}+ \frac{U_{00}\vartheta^{1/k+1}}
{ n^{2(k+1)/(2k+3)}}\big(X_i-\widehat{v}_{i-1}\big)
\\
\widehat{v}^{(j)}_i&=\widehat{v}^{(j)}_{i-1}+
\frac{1}{n}\widehat{v}^{(j+1)}_{i-1}+
\frac{U_{0j}\vartheta^{(j+1)/k+1}} {n^{(2(k+1)-j)/(2k+3)}}
\big(X_i-\widehat{v}_{i-1}\big)
\\
j&=1,\ldots,k-1
\\
\widehat{v}^{(k)}_i&=\widehat{v}^{(k)}_{i-1}
\Big(1-\frac{a_1}{n}\Big)-\frac{1}{n}
\Big(\sum_{\ell=2}^{k-1}a_\ell\widehat{v}^{(k-\ell)}_{i-1}+
a_k\widehat{v}_{i-1}+a_kK\Big)
\\
&\quad + \frac{U_{0k}\vartheta}{n^{(k+2)/(2k+3)}}
\big(X_i-\widehat{v}_{i-1}\big).
\end{aligned}
\end{equation}


\begin{thebibliography}{99}

\bibitem{Torben} Andersen, T., Bollerslev, T., Diebold, F.X. and Labys, P.
{\em Exchange Rate Returns Standardized by Realized Volatility are
(Nearly) Gaussian. }, Multinational Finance Journal , {\bf 4},
2000, pp. 159--179.

\bibitem{BB} Baillie, R.T., T, Bollerslev, {\em Prediction in Dynamic Models with
Time-Dependent Conditional Variances,} Journal of Econometrics,
{\bf 52}, 1992, pp. 91-113.



\bibitem{Black1}  Black, F. {\em The Pricing of Commodity Contracts.
Journal of Financial Economics,} {\bf 9}, 1976, pp. 167-179.

\bibitem{Black}  Black, F. and M. Scholes, {\em The Pricing of Options and Coporate
Liabilities,} Journal of Political Economics, \textbf{81}, 1973,
pp. 637-659.



\bibitem{Bol86} Bollerslev, T., {\em Generalized Autoregressive Conditional
Heteroskedasticity,} Journal of Econometrics, {\bf 31}, pp.
307-327, 1986.

\bibitem{BCK} Bollerslev, T., R.Y. Chou, K.F. Kroner,
{\em ARCH Modeling in Finance: A Review of the Theory and
Empirical Evidence,} Journal of Econometrics, {\bf 52}, pp. 5-59,
1992.




\bibitem{CKL} Chow, P.L., Khasminskii, R. and Liptser, R.  (1997)
 Tracking of signal and its derivatives in the Gaussian White noise, Stoc. Proc. Appl.
 69, 2,  259-273.


\bibitem{Day} Day, T.E. and C.M. Lewis . {\em Forecasting Futures Market
Volatility}, Journal of Derivatives, Winter 1993.

\bibitem{Duan} Duan, J. C.  {\em The GARCH Option Pricing Model.
Mathematical Finance.} {\bf 5} (1), 1995, pp. 13-32.


\bibitem{Eng82} Engle, Robert, {\em Autoregressive Conditional Heteroskedasticity with
Estimates of the Variance of United Kingdom Inflation,}
Econometrica, {\bf 50}, pp. 987-1007, 1982.


\bibitem{GN} Goldenshluger, A. and Nemirovski, A. {\it Adaptive de-noising of signals
satisfying differential inequalities.} IEEE Transactions on
Information Theory, {\bf 43}, 3, pp. 873-889, 1997

\bibitem{GL} Goldentayer, L. and Liptser, R. {\it On-line tracking of a smooth
regression function,}  Statistic Inference for Stochastic
Processes, to apear in 2003.


\bibitem{Ham94} Hamilton, J.D.,
{\em Time Series Analysis,} Princeton University Press, 1994.




\bibitem{IK80} Ibragimov, I. and Khasminskii, R.,  {\em On nonparametric
estimation of regression}, Soviet Math.Dokl.,{\bf 21}, 1980, pp.
810--814.

\bibitem{IK81}  Ibragimov, I.  and  Khasminskii, R.  {\em Statistical estimation:
Asymptotic theory.} Springer Verlag, 1981 (Russian ed.1979).


\bibitem{KhL} Khasminskii, R. and  Liptser, R.  On-line estimation of a
smooth regression function. {\it Theory of Probability and its
Applications}. {\bf 3} (2002).


\bibitem{Spok} Mercurio, Danilo and Spokoiny, Vladimir,
{em Statistical inference for time-inhomogeneous volatility
models.}, {\bf 583}, 2000,
www.wias-berlin.de/publications/preprints/index-2000.html










\bibitem{St82} Stone, C.   {\em Optimal global rates of convergence for
nonparametric regression,} Ann. Statist., {\bf 10}, 1982,
pp.1040--1053.
\end{thebibliography}
\end{document}